\documentclass[12pt]{article}

\usepackage{amsfonts}  
\usepackage{color}

\title{${\cal C}^{1,\beta}$ regularity for Dirichlet problems associated to fully nonlinear degenerate elliptic 
equations.} 

\author{ I. Birindelli, F. Demengel}

\date{}

\catcode`@=11 \@addtoreset{equation}{section} \catcode`@=12

\newtheorem{theo}{Theorem}[section]
\newtheorem{prop}[theo]{Proposition}
\newtheorem{rema}[theo]{Remark}

\newtheorem{cor}[theo]{Corollary}
\newtheorem{lemme}[theo]{Lemma}

\def\R{{\rm I}\!{\rm  R}}
\def\grad{\nabla}

\begin{document}

\maketitle
\section{Introduction}
In a recent paper Imbert and Silvestre \cite{IS} have proved that the solutions of

\begin{equation}\label{eq0}
|\grad u|^\alpha F(D^2 u)=f(x)\ \mbox{in}\ \Omega\subset\R^N
\end{equation}
when $F$ is uniformly elliptic and $f$ is continuous, have first derivative which are H\"older continuous when $\alpha\geq 0$ in the interior of $\Omega$.

Results concerning regularity of solutions have an intrinsic interest which doesn't 
need to be explained. When $\alpha=0$, it is known  (see e.g. Evans \cite{E}, Cabr\'e, Caffarelli \cite{Ca}, 
\cite{CaC1}, \cite{CaC}) that 
$u$ is $C^{1,\beta}$ for some $\beta\in (0,1)$; when $F$ is concave in the
Hessian the solutions are $C^{2,\beta}$.
But for solutions of (\ref{eq0}) with $\alpha>-1$ the question  of the  continuity of the gradient was open and 
it was naturally raised, in \cite{BD9}.  We shall briefly recall why. 
The values
$$\lambda^+=\sup\{\lambda, \exists \phi>0 \ \mbox{in}\ \Omega,  |\grad \phi |^\alpha F(D^2\phi)+\lambda \phi^{1+\alpha} \leq 0\quad\mbox{in}\ \Omega\}$$
$$\lambda^-=\sup\{\lambda, \exists \psi<0 \ \mbox{in}\ \Omega, |\grad \psi |^\alpha F(D^2\psi)+\lambda |\psi|^{\alpha}\psi \geq 0\quad\mbox{in}\ \Omega\}$$
are  generalized principal eigenvalues in the sense that there exists 
a non trivial solution to the Dirichlet problem 
$$ |\grad\phi|^\alpha F(D^2 \phi)+\lambda^{\pm}|\phi|^{\alpha}\phi= 0\quad\mbox{in }\quad \Omega,\ \phi=0 \quad\mbox{on}\quad \partial\Omega,$$
with constant sign.

The main scope in \cite{BD9} was to prove the simplicity of these principal eigenvalues.
The difficulty comes from the fact that the strong 
comparison principle holds only in open subsets of $\Omega$ where the gradient is bounded away from 
zero (in the viscosity sense).  It is well known that the Hopf Lemma guaranties that this is true on
$\partial\Omega$.  So that the continuity of the gradient up to the boundary implies that, in a neighborhood 
of it, the strong comparison principle holds, which is exactly what is needed to prove that the eigenvalues 
are simple.

In that same paper, we proved that if $\alpha\in(-1,0)$ the solutions of the Dirichlet problem
associated to (\ref{eq0}) are indeed $C^{1,\beta}$ and we raised the problem of wether that regularity 
would hold also for $\alpha\geq 0$ i.e. when the operator is degenerate elliptic.

\cite{IS} was a first answer in that direction. Very much inspired by that breakthrough, we wanted to 
complete the work, since, as we have just explained above, in order to use the regularity result in the proof
of the simplicity of the eigenvalues,  it is essential to prove the regularity of the derivative up to the 
boundary.

Recall that $F$ is uniformly elliptic if there exists $\Lambda\geq\lambda>0$ such that for any symmetric 
matrix  $M$
\begin{equation}\label{HF}
{\cal M}^-_{\lambda,\Lambda}(M)\leq F(M)\leq {\cal M}^+_{\lambda,\Lambda}(M)
\end{equation}
here, ${\cal M}^-_{\lambda,\Lambda}(M)=\lambda  tr(M^+)+\Lambda tr(M^-)$ and ${\cal M}^+_{\lambda,\Lambda}(M)=\lambda  tr(M^-)+\Lambda tr(M^+)$. In the rest of the paper we shall drop the indices $\lambda$ and $\Lambda$ of the Pucci operators.

We now  state  the regularity's  result  we prove in this paper.

\begin{theo}\label{th1}
Suppose that $\Omega$ is  a bounded ${\cal C}^2$ domain of $\R^N$ and $\alpha\geq 0$. Suppose that 
$F$ is uniformly elliptic, that $h$ is a continuous function such that
$\left(h(x)-h(y)\right)\cdot (x-y) \leq 0$. Let $f\in{\cal C}(\overline{\Omega})$ and $\varphi\in {\cal C}^{1, \beta_o} (\partial \Omega)$.  For any $u$,  viscosity solution  of 
$$\label{(1)} \left\{ \begin{array}{lc}
   |\nabla u|^\alpha (F(D^2 u) +h(x)\cdot\grad u)= f & {\rm in} \ \Omega\\
    u = \varphi & {\rm on} \ \partial \Omega
    \end{array}\right.$$
there exist $\beta = \beta (  \lambda , \Lambda, |f|_\infty, N, \Omega, |h|_\infty, \beta_o )$  and $C = C(\beta)$ such that 
    $$||u||_{{\cal C}^{1, \beta} (\overline{\Omega})} \leq  C\left( ||\varphi||_{ {\cal C}^{1, \beta_o} (\partial \Omega)}+  |u|_{L^\infty(\Omega)}+ |f|_{L^\infty(\Omega)}^{1\over 1+\alpha}\right).$$
    \end{theo}

For radial solutions, and a more general class of operators, this was proved in 
\cite{BDrad}, with  the optimal H\"older's coefficient $\beta = {1\over 1+\alpha}$.

The novelty with respect to the paper of Imbert and Silvestre is two folded, on one hand we  have added  the  lower 
order term $h(x) \cdot \nabla u |\nabla u|^\alpha $,  and on the other hand we go all the way to the boundary.
The proof follows the scheme of the one in \cite{IS}, but requires new tools and new ideas. In particular
in section 2, we give some a priori Lipschitz and H\"older estimates in the presence of boundary condition 
on one part of the boundary. These are important because the
proof of Theorem \ref{th1} requires that sequence of bounded solutions do converge 
to a solution of a limit equation.
The main tool remains an "improvement of flatness lemma". The reader will see that the presence of the
lower order term and of the boundary term requires some new idea in order to complete the proof.  

For completeness sake, let us now write down the theorem concerned with the simplicity of the
principal eigenvalues.
\begin{theo} Suppose that $\Omega$ is a bounded ${\cal C}^2$ domain,  
such that  $\partial \Omega$ is connected. Then let $\varphi$ and $\psi$ be two  positive eigenfunctions  for 
the principal eigenvalue $\bar \lambda$, i.e.  they are both  solutions of 
            $$\left\{ \begin{array}{lc}
   |\nabla u|^\alpha (F(D^2 u) + h(x)\cdot \nabla u) +  \bar \lambda u ^{1+\alpha} = 0 & {\rm in} \ \Omega\\
    u = 0 & {\rm on} \ \partial \Omega
    \end{array}\right.$$
Then there exists  $t>0$ such that $\varphi = t\psi$.
     
Of course the same result holds for $\underline{\lambda}$ and negative eigenfunctions.
\end{theo}
We will not give the proof of this Theorem, since the  proof given  in \cite{BD9}  in the case 
$\alpha \leq 0$ can be extended to the case $\alpha >0$ as soon as the  eigenfunctions are ${\cal C}^1$ 
near the boundary.

 \section{Local  H\"older and Lipschitz estimates up to the boundary.}

Throughout the paper, the  notation $B_r(x)$ indicates the ball of radius $r$ and center $x$, the center may  be dropped if no ambiguity arise.

It is a classical fact that in order to prove that $u$ is ${\cal C}^{1,\beta}$
at $x_o$, it is enough to prove that  there exists some constant $C$ 
such that, for all $r<1$, there exists $p_r$, such that 
osc$_{B_r(x_o)} (u(x)-p_r\cdot x ) \leq C r^{1+\beta}$. 
   
Furthermore, $u$ is ${\cal C}^{1,\beta}$ in some bounded open set $B$ 
if there exists a constant $C_\beta$ such that for all 
$x\in B$ and $r<1$, there exists $p_{r,x}$ such that 
$$     {\rm osc}_{B_r(x)} (u(y)-p_{r,x}\cdot y ) \leq C_\beta r^{1+\beta}.$$ 
This will be used in  the whole paper. 

We begin by stating the following comparison theorem which will be employed several times.
 
\begin{theo}\cite{BD1}
\label{thcomp}
Suppose that $\Omega$ is a bounded  open set in $\R^N$.  
Let $u$ and $v$ be in ${\cal C}(\overline{\Omega})$  and  respectively solutions of 
$$ |\nabla u|^\alpha (F(D^2 u)+ h(x)\cdot \nabla u) \leq f\ \mbox{ in}\ \Omega$$ 
and 
$$ |\nabla v|^\alpha( F(D^2 v)+ h(x) \cdot \nabla v) \geq  g\ \mbox{ in}\ \Omega$$
with $f$ and $g$ continuous and bounded  such that $f < g$.

If  $u\geq v$ on $\partial \Omega$ then  $u\geq v$ in $\Omega$.
\end{theo}

In order to prove H\"older and Lipschitz estimates we fixe a few notations concerning $\Omega$ and $F$.
We suppose, without loss of generality, that at $0\in \partial \Omega$,  the interior normal is $e_N$. 
By the implicit function theorem, 
there exist a ball $B=B_R(0)$ in $\R^N$,  a ball $B^\prime=B^\prime_R(0)$ of $\R^{N-1}$ and  $a\in {\cal C}^2(B_R^\prime(0))$, 
such that  $a(0)=0$, $\grad a(0)=0$ and, for $y=(y^\prime,y_N)$,
$$\Omega \cap B = \{ y_N> a(y^\prime), y^\prime\in B^\prime \},\  \mbox{and} 
\ \partial \Omega \cap B = \{ y_N=a(y^\prime), y^\prime \in B^\prime\}.$$ 
We shall also act as if $F$ be positively  homogenous of degree 1 i.e. such 
that for any $t>0$, $F(tM)=tF(M)$. 
Observe though that, if this
doesn't hold, when it occurs, we replace  in the computations $F(M)$ by $G_t(M)=t^{-1}F(tM)$; 
this operator satisfies (\ref{HF}) with the same constant than $F$ and the results are unchanged.

In the lemma below we have supposed, for simplicity, that $B$ is the unit ball centered at the origin.

\begin{lemme}\label{lemboundary} Let $\varphi$ be a  H\"older continuous function.
Let  $a\in{\cal C}^2(B^\prime)$ such 
that $a(0)=0$ and $\nabla a (0)=0$. Let $d$ be the distance to the hypersurface $\{ y_N = a(y^\prime)\}$. 

Then,  for all $r<1$ and for all $\gamma <1$, there exists   $\delta_o = \delta_o(|f|_\infty, \gamma, r, \lambda, \Lambda, |h|_\infty, \Omega, \varphi )$,  such that for all $\delta <\delta_o$   
 any $u$,  $|u|_\infty \leq 1$, solution   of 
\begin{equation}\label{eqqq}
\left\{\begin{array}{cc}
       | \nabla u|^\alpha (F(D^2 u)+h(y)\cdot\grad u) = f& {\rm in} \  B\cap \{ y_N >a(y^\prime)\}\\
         u=\varphi &{\rm on} \  \  B\cap \{ y_N = a(y^\prime)\}
       \end{array}\right.
\end{equation}
satisfies $ |u(y^\prime, y_N)-\varphi(y^\prime)| \leq {6\over \delta}  { d(y)\over 1+  d(y)^\gamma}$ in  
$ B_r \cap\{y_N >a(y^\prime)\} $.  
 \end{lemme}
{\em Proof.} 
We write the details of the proof for $\varphi=0$.  The changes to bring in the case where  $\varphi\neq 0$ 
will be given at the end of the proof, the detailed calculation being left to the reader. 

It is  sufficient to consider the set where $ d(y)  < \delta$ since the assumption  
$|u|_\infty\leq 1$ implies the result elsewhere. 

We begin to choose $\delta< \delta_1$,  such that on $d(y) < \delta_1$  the distance 
is ${\cal C}^2$ and  satisfies 
$|D^2 d|\leq C_1$.  We shall also choose later $\delta$ smaller in function of 
$(\lambda, \Lambda, |f|_\infty, |h|_\infty ,  N)$.

Let
$$w(y) = \left\{ \begin{array}{lc}
            {2\over \delta} {d(y) \over 1+ d^\gamma(y)}& \ {\rm for} \  |y  |  < r\\
             {2\over \delta} {d(y)\over 1+ d^\gamma(y)} + {1\over  (1-r)^3} (|y|- r)^3&\ {\rm for} \   | y | > r.
             \end{array}\right.$$
We first remark that 
$$w\geq u\ \mbox{on}\
  \partial (B \cap \{ y_N > a(y^\prime)\}\cap \{d(y)< \delta\}).$$
Indeed, let us observe that on $\{ d(y)=\delta\} $, 
$w\geq  {2\over \delta} {\delta \over 1+\delta^\gamma} \geq 1\geq u$.  
On  $\{|y | = 1\}\cap\{d(y) < \delta\}$, $w \geq {1\over (1-r)^3}(1-r)^3 \geq u $. 
On $B\cap\{ y_N = a(y^\prime)\} $, $w\geq 0 = u$.

We need to choose  $\delta $ small enough that $w$ satisfies
\begin{equation}\label{ww}
| \nabla w|^\alpha ({\cal M}^+ (D^2 w)+h(y)\cdot\nabla w) <-| f|_\infty , \   {\rm in} \  B \cap \{ y_N > a(y^\prime)\}\cap \{d(y)< \delta\}.
\end{equation}
For that aim,  we compute 
  $$\nabla w =\left\{ \begin{array}{lc}
   {2\over \delta}  {1+ (1-\gamma) d^\gamma\over  (1+ d^\gamma)^2} \nabla d & {\rm  when} \ |y|< r\\
  {2\over \delta}  {1+ (1-\gamma) d^\gamma\over (1+ d ^\gamma)^2} \nabla d + {y \over |y| }{ 3\over (1-r)^3}(|y|-r)^2 & {\rm  if} \ |y| > r.
  \end{array}\right. $$
Note that $|\nabla  w|\geq {1\over 4\delta}$ as soon as $\delta\leq\frac{1-r}{12}$. 
By construction $w$ is ${\cal C}^2$ and
$$D^2 w =\left\{ \begin{array}{lc}
     -({\gamma d^{\gamma-1}\over 2\delta} ){(1+\gamma) + (1-\gamma) d^\gamma\over (1+ d^\gamma)^3}\nabla d \otimes  \nabla d +{2\over \delta}  {1+ (1-\gamma) d^\gamma\over  (1+ d^\gamma)^2}   D^2 d  & {\rm if}\  |y| < r\\
     -({\gamma d^{\gamma-1}\over 2\delta} ){(1+\gamma) + (1-\gamma) d^\gamma\over (1+ d^\gamma)^3}\nabla d \otimes  \nabla d +{2\over \delta}  {1+ (1-\gamma) d^\gamma\over  (1+ d^\gamma)^2}   D^2 d +  H(y) & {\rm if}\  |y| > r
     \end{array}\right.$$
where $\|H(y)\|\leq {6\over (1-r)^2}+\frac{3(N-1)}{r(1-r)}$. 
     
We now choose $\delta$ small enough in order that  $\delta < \delta_1$ and 
$$ \lambda (\gamma \delta ^{\gamma-2 } ){(1+\gamma) \over (1+ \delta^\gamma)^3} > 2\Lambda\left({6\over (1-r)^2}+\frac{3(N-1)}{r(1-r)}+  {2C_1\over \delta}\right) + {4|h|_\infty\over \delta}$$ 
and also such that 
$${\lambda\over 2^{2+\alpha}} (\gamma \delta ^{\gamma-(2+\alpha) } ){(1+\gamma) \over (1+ \delta^\gamma)^3} > ||f||_\infty.$$

These choices of $\delta$, and standard computations that use (\ref{HF}) imply that $w$ satisfies (\ref{ww}).
The comparison principle gives that $u\leq w$ in $B \cap \{ y_N > a(y^\prime)\}\cap \{d(y)< \delta\}$. 

Furthermore  the  desired lower bound  on $u$   is easily deduced by considering 
$-w$ in place of $w$ in the previous computations and restricting to $B_r \cap \{ y_N > a(y^\prime)\}$.
This  ends the proof of  Lemma \ref{lemboundary}.

\medskip
When $\varphi\not\equiv 0$, we extend $\varphi$ as a solution  $\psi$ of
$$\left\{\begin{array}{lc}
       {\cal M}^+ (D^2 \psi) = 0& {\rm in} \  B\cap \{ y_N >a(y^\prime)\}\\
         \psi =\varphi & {\rm on} \  B\cap \{ y_N = a(y^\prime)\}
       \end{array}\right.$$
and $\psi$ is ${\cal C}^{1, \beta_o} (B\cap \{ y_N \geq a(y^\prime)\})\cap {\cal C}^2 ( B\cap \{ y_N >a(y^\prime)\})$. 
Furthermore, we can choose $\psi$ such that $|\psi|_\infty\leq |\varphi|_\infty\leq 1$, 
$|\nabla \psi|_\infty \leq c |\nabla \varphi|_\infty$, for some constant which depends on 
$\lambda, \Lambda, N, \Omega$
\cite{CaC}.

We now define
$$w(y) = \left\{ \begin{array}{lc}
            {8\over \delta} {d(y) \over 1+ d^\gamma(y)}+ \psi(y) & \ {\rm for} \  |y|  < r\\
             {8\over\delta} {d(y)\over 1+ d^\gamma(y)} + {1\over  (1-r)^3} (|y|- r)^3+ \psi (y) &\ {\rm for} \   |y| > r.
             \end{array}\right.$$
Similar computations imply that 
$$w\geq u\ \mbox{on}\
  \partial (B \cap \{ y_N > a(y^\prime)\}\cap \{d(y)< \delta\}).$$
Furthermore 
choosing $\delta$  small enough,  we  can ensure that 
$$|\nabla w|^\alpha( {\cal M} (D^2 w)+ h(x) \cdot \nabla w ) < -|f|_\infty.$$
For the lower bound, we replace $w$ by $2\psi-w$. 

\bigskip

Using this estimate together with  an argument due to  Ishii and Lions, \cite{IL},  one finally gets the H\"older  
regularity of the solution, which can be stated as  follows with the same hypothesis on $a$,  $B$  and $f$ 
as above : 

\begin{prop}\label{hold} Let $\varphi$ be a Lipschitz continuous function.
Suppose that $u$ satisfies (\ref{eqqq}).

For all $r<1$,   and for all $\gamma $,  $u$ is $\gamma $ H\"older continuous   on 
$B_r \cap \{ y_N > a(y^\prime)\}$,  with some H\"older's  constant depending on $(r, \lambda, \Lambda, a, N, |f|_\infty ,|h|_\infty,  {\rm Lip} \varphi)$. 
\end{prop}
\begin{rema}\label{remanoboun}
In the absence of boundary conditions,  the solutions are H\"older continuous inside $B_r$ for any 
$r $ such that $B_r \subset \subset B$. 
We do not give the proof which follows the lines in the proof below, it is sufficient to  cancel  in it the 
dependence on $\varphi$. This will be used in the proof of  the interior improvement of flatness lemma with 
additional  lower terms. 
   \end{rema}
{\em Proof  of Proposition \ref{hold}. } We use both some arguments in \cite{BD1} and \cite{IS}. 
Let $1>r^\prime> r$.   Without loss of generality we can suppose that osc $u\leq 1$.
Let $x_o\in B_r \cap \{ y_N > a(y^\prime)\}$ and $\Phi$ defined as 
 
$$\Phi (x, y) =u(x)-u(y)- M |x-y|^\gamma- L |x-x_o|^2- L|y-x_o|^2.$$
The scope is to prove that for $L$ and $M$ independent of $x_o$,  chosen large enough,  
\begin{equation}\label{phi}
\Phi(x, y) \leq 0\ \mbox{on}\ B_{r^\prime} \cap \{ y_N > a(y^\prime)\}.
\end{equation}
This will imply   that $u$ is $\gamma$-H\"older continuous on $B_r \cap \{ y_N > a(y^\prime)\}$ by taking 
$x= x_o$. 

To prove (\ref{phi}), we begin to observe that  it is true on the boundary.
First we suppose that $y_N=a(y^\prime)$.  According to Lemma \ref{lemboundary}  there exists $M_o$  such that  for  $x\in B_{r^\prime} \cap \{ y_N > a(y^\prime)\}$,
$$ |u(x)-\varphi ( x^\prime) |\leq  M_o d(x, \partial \Omega).$$
Then, using 
$|x^\prime - y^\prime | \leq |x-y|$, one has 
 \begin{eqnarray*}
  |u(x^\prime, x_N) - u( y^\prime, a(y^\prime))|&\leq&  |u(x^\prime, x_N) - u( x^\prime, a(x^\prime))|+ |u(x^\prime, a(x^\prime)) - u( y^\prime, a(y^\prime))|\\
  &\leq& M_o d(x,\partial\Omega)+{\rm Lip}_\varphi|x^\prime-y^\prime|\\
  &\leq& M_o |x-(y^\prime, a(y^\prime))|+{\rm Lip}_\varphi |x-(y^\prime, a(y^\prime))|.
  \end{eqnarray*}
So, if $M$ is chosen greater than $ M_o+ {\rm Lip} \varphi$,  then we have obtained that 
$\Phi\leq 0$ on 
$B_{r^\prime} \cap \{ y_N = a(y^\prime)\}$.

On the rest of the boundary it is enough to choose $L > {4\over (r^\prime- r)^2}$ and to recall that the oscillation of $u$ is bounded by 1. 

\noindent In the sequel we will choose $M$ large in order that  
${L\over M} = o(1)$.  

Suppose by contradiction that $ \Phi(x, y)>0$
for some $(x,y)\in B_{r^\prime }   \cap \{ y_N > a(y^\prime)\}$.            
           Then there exists $(\bar x, \bar y)$ such that 
$$\Phi(\bar x, \bar y)=\sup_{\overline{B_{r^\prime}}} (\Phi(x,y))>0.$$
 Clearly $\bar x\neq \bar y$. Furthermore the hypothesis on $L$ forces $\bar x$ and $\bar y$ to be in $B_{r^\prime+ r\over 2} \cap \{ y_N > a(y^\prime)\}$. 
Then,  for all $\epsilon >0$ small depending on the norm of $Q := D^2 (M |x-y|^\gamma)$,  using Ishii's Lemma \cite{I},  there exist $X$ and $Y$ such that
$$( \gamma M(\bar x-\bar y)|\bar x-\bar y|^{\gamma-2}+ 2L(\bar x-x_o)
, X)\in J^{2,+ }u^\star(\bar x)$$
$$(\gamma M(\bar x-\bar y)|\bar x-\bar y|^{\gamma-2}- 2L(y-x_o)
, -Y)\in J^{2,-}u_\star(\bar y)$$
with 
$$\left(\begin{array}{cc}
X&0\\
0&Y
\end{array}\right) \leq \left(\begin{array}{cc}
Q &-Q\\
-Q&Q
\end{array}\right) +(2L +\epsilon) \left(\begin{array}{cc}
I&0\\
0&I\end{array} \right). $$

In the sequel, since we assumed that ${L\over  M} = o(1)$, one also has 
${L+\epsilon \over M}= o(1)$ and then we drop $\epsilon $ for simplicity.

 Let us denote  $q_x=  \gamma M(\bar x-\bar y)|\bar x-\bar y|^{\gamma-2}+ 2L(\bar x-x_o)$, and  $q_y =  \gamma M(\bar x-\bar y)|\bar x-\bar y|^{\gamma-2}-2L(\bar y-x_o)$. 
By the choice of $M$ in function of $L$, 
$$ 2\gamma M |\bar x-\bar y|^{\gamma-1} \geq ( |q_x|,\ |q_y|) \geq {1\over 2 } \gamma M |\bar x-\bar y|^{\gamma-1},$$  
and since 
$|q_x-q_y|\leq  4L$,  by the mean value's theorem and using a constant $\kappa  < \alpha$ 
if $\alpha <1$ and $\kappa=1$ if $\alpha >1$: 
     \begin{eqnarray*}
     ||q_x|^\alpha -|q_y|^\alpha |&\leq& \alpha |q_x-q_y|^{\kappa}  2^{|\alpha-1|}   |\gamma M |\bar x-\bar y|^{\gamma-1} |^{1-\kappa } |\gamma M |\bar x-\bar y|^{\gamma-1} |^{\alpha-1} \\
     &\leq &C (M\gamma|\bar x-\bar y|^{\gamma-1})^{\alpha-\kappa} \left({L\over M}\right)^\kappa = o((M|\bar x-\bar y|^{(\gamma-1)})^{(\alpha-\kappa)}).
     \end{eqnarray*}
We now treat the  terms concerning the second order derivative.
The previous inequalities can also be written as 
$$\left(\begin{array}{cc}
X-2LI &0\\
0&Y-2L I
\end{array}\right) \leq \left(\begin{array}{cc}
Q &-Q\\
-Q&Q
\end{array}\right).$$
 We prove in what follows that,  for some constant which can vary from one line to another,
 $L = o(|tr(X+Y)|) $,  and that  there exist constant $c$ and $C$ such that 
 $|X|, |Y| \leq C |tr(X+ Y)|  $  and 
 $$cM |\bar x-\bar y|^{\gamma-2} \leq |tr(X+Y)|\leq CM  |\bar x-\bar y|^{\gamma-2}.$$

Indeed, let 
$$ P : = {(\bar x-\bar y\otimes \bar x-\bar y)\over |\bar x-\bar y|^2}\leq I.$$

Using $-(X+Y)\geq 0$ and $(I-P)\geq 0$ and the properties of the
symmetric matrices one has 
$$tr(X+Y-4L)\leq tr(P(X+Y-4L)).$$ 
Remarking in addition that
$X+Y-4L\leq 4Q$, one  sees that
$tr(X+Y-4L)\leq 4tr(PQ)$. But $tr (PQ)=\gamma M(\gamma-1)|
\bar x-\bar y|^{\gamma-2}<0$, hence

\begin{equation}\label{elip}
|tr (X+Y-4L)|\geq 4\gamma M(1-\gamma)| \bar x-\bar y|^{\gamma-2}.
\end{equation}
Furthermore  by Lemma III.1 of \cite{IL} there exists a universal constant $C$ such that 
\begin{eqnarray*}
|X|, |Y|, |X-2L|, |Y-2L|&\leq& C (|tr(X+Y-4L)|+ |Q|^{1\over 2} |tr(X+Y-4L)|^{1\over
2})\\
&\leq& C|tr(X+Y-4L)|\\
&\leq & C|tr(X+Y)|
\end{eqnarray*}
since $|Q|$ and $|tr(X+Y-4L)|$ are of the same order, and ${L\over M}  = o(1)$. 
This will yield the required  estimates.

For some  positive constants  $c_2$, $c_3$,  since $u$ is both a sub-  and a 
super-  solution of  (\ref{eqqq}), using the uniform ellipticity of $F$ and the  assumptions on  $h$:
\begin{eqnarray*}
f(\bar x) &\leq &  |q_x|^\alpha( F(X)+ h(\bar x)\cdot q_x)\\
                &\leq& |q_y|^\alpha (F(X)+ h(\bar x) \cdot q_x) \\
                &&+ o(M\gamma |\bar x-\bar y|^{\gamma-1})^{\alpha-\kappa} (\Lambda |X|+|h|_\infty (\gamma M|\bar x-\bar y|^{\gamma-1}+2 L))\\
                &\leq &   |q_y|^\alpha( F(Y)  + h(\bar y) \cdot q_x+ 4|h|_\infty L) +\\
                &&+ |q_y|^\alpha tr(X+Y) +   o( M^{1+\alpha}|\bar x-\bar y|^{(\gamma-1 )\alpha+ \gamma-2})\\
                               &\leq & M^\alpha c_2 |\bar x-\bar y|^{(\gamma-1)\alpha}  tr(X+ Y)+ o( M^{1+\alpha} |\bar x-\bar y|^{(\gamma-1 )\alpha+ \gamma-2})+ f(\bar y)\\\
               &\leq &-c_3 M^{1+\alpha} |\bar x-\bar y|^{(\gamma-1) \alpha + \gamma -2} +  f(\bar y).\\
 \end{eqnarray*}
This is clearly false  as soon as $M$  is large enough and it ends the proof. 
              
\bigskip

Compactness near the boundary is a natural consequence of Proposition \ref{hold}. 
\begin{cor}\label{compac1}
Suppose that $(u_n)$  is a bounded sequence of continuous functions which satisfy  
$$\left\{ \begin{array}{lc}
           |\nabla u_n|^\alpha (F( D^2 u_n) +h(y)\cdot\grad u_n)= f_n  & {\rm in} \ B\cap  \{ y_N >a(y^\prime)\}\\
            u_n = \varphi& {\rm on } \ B \cap \{ y_N = a(y^\prime)\}
 \end{array} \right.$$
and suppose that $(f_n)$ converges simply to some continuous function $f$. Then  for all $r<1$,  
 one can extract from $(u_n)$  a subsequence which converges uniformly,  on  $\overline{B_r  \cap \{ y_N >  a(y^\prime)\}}$,  towards a solution of 
              $$\left\{ \begin{array}{lc}
           |\nabla u|^\alpha (F( D^2 u)+h(y)\cdot\grad u) = f  & {\rm in} \ B\cap \{ y_N >a(y^\prime)\}\\
            u = \varphi& {\rm on } \ B \cap \{ y_N = a(y^\prime)\}.
            \end{array} \right.$$
                        \end{cor}
                          \begin{rema}\label{compnoboun}
In the absence of boundary conditions,  the  analogous result holds, in the sense that one can extract from 
$(u_n)_n$ a subsequence which converges uniformly on every $B_r\subset \subset B$  to  $u$ a solution of        
$$ |\nabla u|^\alpha (F( D^2 u)+h(y)\cdot\grad u) = f  \ {\rm in} \ B.$$
   \end{rema}
 When we shall treat,  in the improvement of flatness lemma up to the boundary, the case where the
 boundary is locally straight, we shall need the following 
Lipschitz estimate's   near the boundary  for some different but related equation.  
 \begin{prop}\label{liplargep}
         (Lipschitz estimates for large $p's$)
 Let $\varphi$ be a Lipschitz continuous function.   Assume that $u$ solves 
          
 $$\left\{ \begin{array}{lc}
 |pe_N+\nabla u|^\alpha (F(D^2 u)+h(y)\cdot \grad u) = f & \ {\rm in}\  B(x) \cap \{ y_N >0\} \\
 u=\varphi & {\rm on} \   \{ y_N=0\}\cap B_1(x) 
 \end{array}\right.$$  
  with $|u|_{L^\infty(B_1(x) \cap \{ y_N >0\} )} \leq 1 $ and $||f||_{L^\infty (B(x) \cap \{ y_N >0\}  )} \leq \epsilon_o <1$. Then, for all $r<1$,
  there exists $b_o$ depending on $(\lambda, \Lambda, N, \alpha, r, \epsilon_o,  {\rm Lip } \varphi)$, 
  such that  if $|p|  > {1\over  b_o}$,  $u$ is Lipschitz continuous in $B_r(x)\cap\{ y_N >0\} $ with some Lipschitz constant   depending on 
  $(\lambda, \Lambda, N, \alpha, r, \epsilon_o, {\rm Lip } \varphi)$. 
    \end{prop}
\begin{rema}\label{remanoboun1}
In the absence of boundary conditions, the solutions are Lipschitz with Lipschitz constant 
independent of $p$,  inside $B_r$,  for any $r $ such that $B_r \subset \subset B$. 
               
This will be used in the proof of  the interior improvement of flatness lemma with   lower order terms. 
\end{rema}
This Proposition  is a consequence of  the following
 \begin{lemme}\label{lemboundary2} For all $\gamma <1$,   for all $r <1,$   there exists  $\delta = \delta (|f|_\infty, \lambda , \Lambda, r) $,   such that  for $b< {\delta \over 4}$,  any solution $u$ of 
$$\left\{\begin{array}{lc}
       |e_N+ b \nabla u|^\alpha (F (D^2 u)+h(y)\cdot\grad u) = f& {\rm in} \  B(x) \cap \{ y_N >0\}\\
         u=\varphi & {\rm on} \  B(x)\cap\{ y_N=0\}.
       \end{array}\right.$$
such that osc$(u) \leq 1$ satisfies $ |u(y^\prime, y_N)-\varphi(y^\prime)| \leq {2\over \delta} {y_N\over 1+ y_N^\gamma}$ in 
$ B_r \cap\{y_N >0\} $.
\end{lemme}
{\em Proof of Lemma \ref{lemboundary2}.} 
Suppose for simplicity that $\varphi=0$. If $b=0$ the result   is known by properties of 
solutions of  $F(D^2 u) + h(x) \cdot \nabla u= f$ which are zero on the boundary. 
So we assume in what follows that $b\neq 0$.

We act as in Lemma \ref{lemboundary},  where the distance is replaced by $y_N$ so we consider
$$w(y) = \left\{ \begin{array}{lc}
 {2\over \delta} {y_N\over 1+ y_N^\gamma}& \ {\rm for}\ y_N < \delta, |y^\prime  |< r\\
{2\over \delta} {y_N\over 1+ y_N^\gamma} + {1\over  (1-r)^3} (|y^\prime |- r)^3&\ {\rm for} \ y_N < \delta,    
|y^\prime| > r.
\end{array}\right.$$
Similarly  to the proof of Lemma  \ref{lemboundary}, it is sufficient to consider the set 
where $y_N< \delta$,  since 
the assumption  $|u|_\infty\leq 1$ implies the result elsewhere. 
Furthermore  we only prove that $u\leq w$, the desired lower bound can be obtained by considering $-w$ in place of $w$.
          
In order for $w$ to satisfy
$$|e_N+ b \nabla w|^\alpha ({\cal M}^+ (D^2 w)+h(y)\cdot\grad w) \leq -| f|_\infty , \   {\rm in} \  B,$$
it is sufficient to choose   $\delta $ such that 
$$\left({1\over 2}\right)^\alpha \lambda \gamma \delta^{\gamma-2} (1-\gamma) {1\over  (1+\delta^\gamma)^3}> |f|_\infty +2\Lambda\left({6\over (1-r)^2}+\frac{3(N-1)}{r(1-r)}+ {4|h|_\infty\over \delta}\right),$$
 $b < {\delta \over 4}$, and recall that $|\nabla w |\leq {2\over \delta}$.
 Furthermore $w\geq u$ on $ \partial(B^+\cap \{0<x_N < \delta\})$.  

The comparison principle   in Theorem \ref{thcomp}, between the functions $x \mapsto  x_N + b w(x)$ and 
$x\mapsto x_N+ b u(x)$, as well as $b\neq 0$,  implies that $u\leq w$ in $B \cap \{ y_N >0\} $. 
Finally the desired estimate is obtained in $\{|y^\prime | < r,\ y_N >0\}$.

In the case $\varphi\neq 0$, we take the function $w$ as in the proof of Lemma \ref{lemboundary}, with 
$d$ replaced by $y_N$. Requiring sufficient restriction on the smallness of $\delta$  give the result. 

\bigskip
We are now ready to give the 

\noindent{\em Proof of Proposition \ref{liplargep}}. We act as in \cite{IS}.  We rewrite the equation as 
$$|e_N+ b Du|^\alpha (F(D^2 u)+h(y)\cdot\grad u) = \tilde f$$ 
with $b ={1\over p}$ and $\tilde f  =   |p|^{-\alpha} f$. 

Choose first $\delta$ small enough in order that 

\begin{eqnarray*}
\left({1\over 2}\right)^\alpha   \lambda (\gamma \delta ^{\gamma-2} ){(1+\gamma) \over  2(1+ \delta^\gamma)^3} &>& |f|_\infty +2\Lambda({6\over (1-r)^2}+\frac{3(N-1)}{r(1-r)}+ {4|h|_\infty\over \delta})\\
&> &|\tilde f|_\infty+2\Lambda({6\over (1-r)^2}+\frac{3(N-1)}{r(1-r)}+ {4|h|_\infty\over \delta}).
\end{eqnarray*}
 and  such that $b < {\delta\over 4}$. 
      Let us note that this implies that $b$  is small enough  depending   on  $(\lambda, \Lambda, N, \alpha, r, \epsilon_o, |h|_\infty, {\rm Lip} \varphi)$. 
      
      \noindent Let $r<r^\prime < 1$ and let  $x_o\in B_{r^\prime }(x) \cap \{ y_N >0\}$,  $L_2= {4\over (r^\prime-r)^2}$,
$$\psi(z,y) = u(z)-u(y)-L_1\omega (|z-y|) -L_2|z-x_o|^2 -L_2 |y-x_o|^2$$
where $\omega(s) = s-\omega_o s^{3\over 2}$ if $s \leq s_o =\left( {2\over 3 \omega_o}\right)^2$ and $\omega (s) = \omega(s_o)$ if $s\geq s_o$.  We also require $L_1> {2\over \delta}+{\rm Lip}_\varphi$.

If we prove that $\psi(z,y)\leq 0$ in $B_{r^\prime} $ , since $L_1$ is independent of $x_o$, 
by choosing $z= x_o$ one gets
$$u(z)-u(y) \leq L_1 |z-y|+ L_2 |z-y|^2$$ which implies the desired result  when $z\in B_{r}(x) $  or $y\in  B_{r}(x)$. 
         
We begin to observe that if the supremum is achieved on $\overline{B_r(x) }$ and if $(\bar x, \bar y)$ is 
a point where the supremum is achieved, then, with our choice of $L_1$, neither $\bar x$ nor $\bar y$ can 
belong to the part $\{z_N = 0\}$  according to  Lemma \ref{lemboundary2}. The rest of the proof is as in 
\cite{IS},  see also the proof of Proposition \ref{hold}.     

\bigskip

As a corollary of this   Lemma one has the following compactness result

\begin{cor}\label{compact2} Let $\varphi$ be a Lipschitz continuous function.
Suppose that $(u_n)$  is a sequence of continuous viscosity solutions of 
          $$\left\{ \begin{array}{lc}
           |e_N+ b_n\nabla u_n|^\alpha\left( F( D^2 u_n)+ h \cdot \nabla u_n\right) = f_n  & {\rm in} \ B,\\
            u_n = \varphi& {\rm on } \ B \cap \{ y_N = a(y^\prime)\}.
            \end{array} \right.$$
           where $b_n\leq b_o$,  $b_o$ is given above in Proposition  \ref{liplargep}. 
Suppose that $f_n$ converges simply to some function $f$  in  $\Omega$.  Then  for all $r<1$,   one can extract from $(u_n, b_n)$  a subsequence which 
converges uniformly on  $\overline{B_r  \cap \{ y_N >  a(y^\prime)\}}\times \R$, and the limit $(u, \bar b)$ 
satisfies  
 $$\left\{ \begin{array}{lc}
           |e_N+ \bar b \nabla u|^\alpha\left( F( D^2 u)+ h\cdot  \nabla u\right)= f  & {\rm in} \ B\\
            u = \varphi& {\rm on } \ B \cap \{ y_N = a(y^\prime)\}.
            \end{array} \right.$$
\end{cor}
   \begin{rema}\label{compnoboun1}
In the absence of boundary conditions, the conclusion is that  the sequence $(u_n)$ contains a 
subsequence which converges locally uniformly and up to a constant toward  a solution of 
$$|e_N+ \bar b \nabla u|^\alpha\left( F( D^2 u)+ h\cdot  \nabla u\right)= f  \  {\rm in} \ B.
$$
 \end{rema}

\section{ Proof of Theorem \ref{th1}.}
In fact Theorem \ref{th1} is an immediate consequence of the following local result up to the boundary 
together with  some argument of finite covering:
     
 \begin{theo}\label{th2} Suppose that $F$, $h$ and $f$ are as in Theorem \ref{th1}.
       Let $B$ be an open set in $\R^N$ and let  $a$ be a ${\cal C}^2$ function defined on  $B\cap \R^{N-1}\times \{0\}$ with $a(0) =0$, $\nabla a (0)=0$. There exists $\beta$ such that for any $u$ solution of 
       $$\left\{ \begin{array}{lc}
   |\nabla u|^\alpha (F(D^2 u) +h(x)\cdot\grad u)= f & {\rm in} \  B \cap \{ x_N > a(x^\prime)\}\\
    u = \varphi & {\rm on} \  B \cap \{ x_N = a(x^\prime)\}, 
    \end{array}\right.$$
    $u$ is ${\cal C}^{1, \beta} (\tilde B \cap \{ x_N > a(x^\prime)\})$ for any $\tilde B\subset\subset  B$.
\end{theo}

Theorem \ref{th2} is proved via the following  two  "improvement of flatness" lemma and their  consequences.
 \begin{lemme}\label{flat}
There exist $\epsilon_o\in ]0,1[$ and $\rho\in ]0,1[$ depending on $(\alpha, |h |_\infty, \lambda, \Lambda, N)$ 
such that for any $p\in \R^N$ and for any viscosity solution $u$  of 
 $$|p+\nabla u|^\alpha \left(F(D^2 u) + h(y) \cdot (\nabla u+p)\right)  = f\ \mbox{in}\ B_1$$
such that $osc_{B_1} u\leq 1$ and $\|f\|_{L^\infty (B_1)} \leq \epsilon_o$, there exists $p^\prime \in \R^N$ such that 
$$osc_{B_\rho} (u-p^\prime \cdot x) \leq {1\over 2} \rho. $$
          \end{lemme}
and
\begin{lemme}\label{lem1a} 
For any $a\in {\cal C}^2$ ,  such that $a(0)=0$ and $\nabla a (0) = 0$,   there exist $\epsilon_o>0$ and $\rho$  which depend on $(\alpha, \lambda, \Lambda, N, |D^2 a|_\infty, |h|_\infty,|\varphi|_{{\cal C}^{1,\beta}}) $ such that for any $p\in \R^N$ and $u$ a viscosity solution of  
$$\left\{ \begin{array}{lc}
 |p+\nabla u|^\alpha (F(D^2 u)+h(y)\cdot(\grad u+p)) = f& \ {\rm in}\   B \cap \{ y_N > a(y^\prime)\}\\
 u+p\cdot y =\varphi & {\rm on} \  \{ y_N=a(y^\prime)\}\cap  B.
 \end{array}\right.$$
  Then for all $x\in B$ such that $B_1(x) \subset B$, 
  $osc_{B_1(x) \cap \{y_N >a(y^\prime)\}}  u \leq 1$, and $|f|_{L^\infty(B_1(x)\cap \{y_N >a(y^\prime)\})}\leq \epsilon_o$,     there exists $q_{x, \rho} \in \R^N$ such that 
  
$$osc_{B_\rho(x)\cap \{y_N >a(y^\prime)\}}  (u(y)-q_{x,\rho}\cdot y) \leq {\rho\over 2}.$$
\end{lemme}

Suppose that these Lemmata  have been proved and let us derive the following one.

  \begin{lemme}\label{lem2a}
Suppose that  $\rho $  and $\epsilon_o \in [0,1]$  are as in lemma \ref{lem1a} 
 and suppose that $u$ is a viscosity solution of 
\begin{equation}\label{eq2a}
\left\{ \begin{array}{lc}
  |\nabla u|^\alpha (F(D^2 u)+h(y)\cdot\grad u) = f& {\rm in } \  B_1(x) \cap\{ y_N >a(y^\prime)\} \\
   u=\varphi& {\rm on } \ B_1(x) \cap \{ y_N=a(y^\prime)\} 
   \end{array}\right.
\end{equation}
 with {\rm osc}\,$u\leq 1$ and $||f||_\infty \leq \epsilon_o$,  then, there exists $\beta \in ]0,1[$, such that for all $k$ and for all $x\in \Omega 
 $ there exists $p_k\in \R^N$ such that  
  $${\rm osc}_{B_{r^k}(x)\cap \{ y_N > a(y^\prime)\} } (u(y)-p_k\cdot y) \leq r_k^{1+\alpha} \equiv {\rho}^{k(1+\alpha)}.$$
 \end{lemme}
 \begin{rema}
  Of course, the interior regularity becomes : 
 
Suppose that  $\rho $  and $\epsilon_o \in [0,1]$  are as in lemma \ref{flat} 
 and suppose that $u$ is a viscosity solution of 
     $$
  |\nabla u|^\alpha (F(D^2 u)+h(\cdot )\cdot\grad u) = f\  {\rm in } \  B_1$$
   with osc$u\leq 1$ and $||f||_\infty \leq \epsilon_o$, then  there exists $\beta \in ]0,1[$ such that for all $k$ 
  there exists $p_k\in \R^N$ such that  
  $${\rm osc}_{B_{r^k}} (u(x)-p_k\cdot x) \leq r_k^{1+\alpha} \equiv {\rho}^{k(1+\alpha)}.$$
   We will not give the proof,  the changes to bring  to the   proof below being obvious. 
 \end{rema}
{\em Proof of Lemma \ref{lem2a}. }
As in \cite{IS} we use a recursive argument.

We first remark that  one can assume that     
$\varphi (x^\prime) =\partial_i \varphi (x^\prime) = 0$ for $i=1,\cdots, N-1$. 

Indeed, let  $u$ be a solution of (\ref{eq2a}). 
Let $v(y) = u(y)-u(x)-\grad\varphi(x^\prime)\cdot(y^\prime-x^\prime)$ which satisfies
$$\left\{ \begin{array}{lc}
 |\nabla v+ q|^\alpha (F(D^2 v)+h(y)\cdot(\grad v+q)) = f& \ {\rm in} \ B_1(x) \cap \{y_N > a(y^\prime)\}\\
v(y^\prime, a(y^\prime))= \phi(y^\prime) & {\rm on} \ B_1(x) \cap \{ y_N = a(y^\prime)\}
\end{array}\right.$$
where $q=p+(\grad\varphi(x^\prime),0)$ and 
$\phi(y^\prime)=\varphi(y^\prime)-\varphi(x^\prime)-\grad\varphi(x^\prime)\cdot(y^\prime-x^\prime)$ which
satisfies $\phi (x^\prime) =\partial_i \phi (x^\prime) = 0$ for $i=1,\cdots, N-1$. 

So the result obtained for $v$ would transfer to $u$ replacing $p$ with $q$.

We can start. For $k=0$, taking $p_o=0$ yields the desired inequality. 
Suppose that $p_k$ has been constructed. 

Choose $\beta$ small enough in order that $\rho^\beta > {1\over 2}$. 

With the above assumption on $\varphi$, let
$\varphi_k (y^\prime) = r_k^{-1-\beta} \varphi ( x^\prime + r_k (y^\prime -x^\prime))$, 
which satisfies, for $\beta < \beta_o$, 
$|\varphi_k|_{{\cal C}^{1, \beta}(B_1(x^\prime))} \leq |\varphi |_{{\cal C}^{1, \beta}} $. 

We consider 
$$u_k (y) = r_k^{-1-\beta} \left(u(r_k (y-x)+ x) -p_k\cdot (r_k (y-x)+x)\right).$$  
$u_k$ is well defined on 
$B_1(x)  \cap \{ y_N > a_k(y^\prime) \}$, where $a_k(y^\prime)=x_N (1-{1\over r_k}) + {a(r_k(y^\prime -x^\prime)+ x^\prime) \over r_k}$.
Note that on 
$B_1(x)  \cap\{ y_N > a_k(y^\prime) \}$, 
$|x_N | < r_k(1+ |D^2 a|_\infty)$. 

It is immediate to see that $u_k$ is a solution of
$$\left\{ \begin{array}{lc}
|p_k r_k^{-\beta} +\nabla u_k|^\alpha (F(D^2 u_k)+h_k\cdot (p_k r_k^{-\beta} +\nabla u_k)) = f_k &\ {\rm in}\   B_1(x) \cap \{ y_N > a_k(y^\prime)\}\\
u_k+p_kr_k^{-\beta} \cdot y =p_k r_k^{-\beta} (r_k-1) \cdot x+ \varphi_k (y^\prime)   & {\rm on} \  
B_1(x)\cap \{ y_N=a_k(y^\prime)\}
 \end{array}\right.$$
with $f_k (y) = r_k^{1-\beta (1+\alpha)} f(r_k (y-x)+x)$ and $h_k(y)=r_kh(r_k (y-x)+x)$.

We  now prove that 
$a_k$ satisfies 
$\|a_k\|_{{\cal C}^2 } \leq C$ for some constant  $C$ which does not depend on $k$. Indeed,
 in $B_1(x)  \cap \{ y_N > a_k(y^\prime)\}$, $r_k  \geq{ |x_N|\over 1+ |D^2 a|_\infty} $,  which implies that
$$|a_k(y^\prime )|\leq   |x_N|+(1+ ||D^2 a||_\infty )+ {r_k^2 ||D^2 a||_\infty \over r_k} \leq {\rm diam } \Omega+2(1 +  ||D^2 a||_\infty)$$
and $ |\nabla a_k| =|\nabla a|$, finally 
$||D^2 a_k||=r_k ||D^2 a||_\infty\leq ||D^2 a||_\infty$. 
           
Furthermore, as long as $\beta < {1\over 1+\alpha}$, 
$${\rm osc}_{B_1(x) \cap \{y_N >a_k(y^\prime)\}}  u_k\leq 1,\ 
 ||f_k||_{L^\infty (B_1(x)  \cap \{ y_N >  a_k(y^\prime)} \leq \epsilon_o\ \mbox{and}\ \|h_k\|_\infty\leq \|h\|_\infty.$$

Hence, using Lemma \ref{lem1a}   with obvious changes,  there 
exists $q_{k+1}  \in \R^N$ such that 
  
  $$\mbox{osc}_{B_\rho(x)\cap \{y_N >a_k(y^\prime)\}}  (u_k(y)-q_{k+1} \cdot y) \leq {\rho\over 2}.$$
Defining $p_{k+1} = p_k + q_{k+1} r_k^\beta$, with the assumptions on $\beta$ and $\rho$, one gets: 
   
$$
{\rm osc}_{ B_{r_{k+1}}(x) \cap \{ y_N > a(y^\prime)\}} \left(u(y)-p_{k+1}\cdot y \right)
   \leq {\rho\over 2} r_k^{1+\beta} \leq r_{k+1}^{1+\beta},
$$
since the oscillation is invariant by translation.
This is the desired conclusion. 
  
    \bigskip
    
There remains to prove the flatness lemmata.  We start by the interior case with lower order terms. 
   
{\em Proof of Lemma \ref{flat}.}   Suppose by contradiction that there exist a sequence of functions $(f_n)_n$ 
whose norm go to zero, a sequence of $(p_n)_n\in\R^N$ and a sequence of functions 
$(u_n)_n$ with osc $u_n\leq 1$,
solutions of 
\begin{equation}\label{pn}
|p_n+\nabla u_n|^\alpha \left(F(D^2 u_n) + h(y) \cdot (\nabla u_n+p_n)\right)  = f_n, 
\end{equation}
such that,  for all $q\in \R^N$, 
\begin{equation}\label{st} {\rm osc}_{B_\rho}(u_n(y)-q\cdot y)> {\rho\over 2}.
\end{equation}

Let us suppose first that $(p_n)_n$ is bounded so,  then up to subsequence, it converges to $p_\infty$. 
Considering $v_n(y) = u_n(y)+ p_n\cdot y$ and using the 
compactness Remark \ref{compnoboun}, we can  extract  form $(v_n)_n$ a subsequence converging to  a limit $v_\infty$, which satisfies    
$$|\nabla v_\infty |^\alpha (F(D^2 v_\infty )+ h(x) \cdot \nabla v_\infty ) =0.$$  
Remark  next that the solutions of such  an equations are 
solutions of 
$$F(D^2 v_\infty)+ h(\cdot ) \cdot \nabla v_\infty =0.$$ 
as it is the case  for $h=0$ (see \cite{IS}). 
But, passing to the limit in  (\ref{st}) gives  
that osc$_{B_\rho}(v_\infty-(q-p_\infty) \cdot x)> {\rho\over 2}$. This contradicts the regularity results known 
for solutions of this equation (see \cite{T}) and it ends the case where the sequence $(p_n)_n$ is bounded.
     
In the case where $(p_n)$ is unbounded, take a subsequence 
such that ${p_n\over |p_n|} $ converges to some $p_\infty$.  

{\bf Claim:}  There exist $q_\infty\in\R^N$ and a subsequence $\sigma(n)$ such that
$$\lim_{n\rightarrow \infty} h(y)\cdot p_{\sigma(n)}=h(y)\cdot q_\infty$$
uniformly in $B_r$ for any $r<1$. 

 We postpone the proof of that claim and end the proof of  Lemma \ref{flat}.

We now divide the equation (\ref{pn}) by $|p_n|^\alpha$ and get, with $e_n = {p_n\over |p_n|}$  and 
$a_n = {1\over |p_n|}$, 
        $$|a_n\nabla u_n+ e_n|^\alpha  \left(F(D^2 u_n) + h(y) \cdot ( \nabla u_n+ p_n) \right) = {f_n \over |p_n|^\alpha}.$$
        
Using   Remark \ref{compnoboun1}      
and the claim, a subsequence of $u_{\sigma(n)}$ converges to $u_\infty$ a solution of
the limit equation
$$ F(D^2 u_\infty ) + h(y) \cdot (\nabla u_\infty+ q_\infty) = 0.$$
On
the other hand,
osc $(u_\infty -q^\prime \cdot x) > {1\over 2} \rho$, 
which contradicts the regularity of solutions of such equations.

{\bf  Proof of the Claim.}  Let $V\subset \R^N$ be the space generated by $h(B_r)$. 
Let $q_n=\Pi_V p_n$ be the projection of $p_n$ on $V$, hence $h(x)\cdot p_n=h(x)\cdot q_n$.

Suppose by contradiction that the sequence $(q_n)_n$ goes to infinity in norm.
There exists a subsequence $q_{\sigma (n)}$ such that
$\frac{q_{\sigma (n)}}{|q_{\sigma (n)}|}\rightarrow  \bar q$, for some $\bar q$ such that $|\bar q|=1$
and  $\bar q\in V$, furthermore we have that  $\frac{p_{\sigma (n)}}{|p_{\sigma (n)}|}\rightarrow  \bar p$ such that $|\bar p|=1$ and $\frac{|q_{\sigma (n)}|}{|p_{\sigma (n)}|}\rightarrow  \tilde p\in\R$.

We divide the equation by $|p_n|^{\alpha}|q_n|$ and observe that
the functions $v_n=\frac{u_n}{|q_n|}$ satisfy

$$ \left|{p_n\over |p_n|}+ {|q_n|\over |p_n|}\nabla v_n\right|^\alpha \left(F(D^2 v_n) + 
{h(y) \cdot q_n\over |q_n| }+h(y)\cdot \grad v_n\right) = {f_n\over |p_n|^{\alpha}|q_n|}.$$
And the sequence $v_n$ converges to zero.

Using  the compactness of $(v_{\sigma(n)})_n$, see Remark \ref{compnoboun1}, one 
gets that it  converges to a solution of
$$(\bar p+\tilde p\nabla v_\infty)^\alpha (F(D^2 v_\infty) + h(y)\cdot\bar q+h(y) \cdot \nabla v_\infty) =0$$ 
but, since $v_\infty=0$ we get
$$ h(y)\cdot \bar q= 0\ \mbox{ for all }\ y\in B_r.$$
This means that $\bar q$ is both in $V$ and $V^\perp$, hence  $\bar q=0$,  
which is a contradiction and the sequence $(q_n)_n$ is bounded.
This gives the claim, indeed, up to a subsequence $q_n$  converges to $q_\infty$ and

$$\lim_{n\rightarrow \infty} h(y)\cdot p_{\sigma(n)}=\lim_{n\rightarrow \infty} h(y)\cdot q_{\sigma(n)}=h(y)\cdot q_\infty.$$
This ends the proof.
\bigskip   

We now want to prove  Lemma \ref{lem1a}. In the case of a non straight boundary it  requires  the 
following technical proposition, whose 
proof is  postponed until the end of the section.
   \begin{prop}\label{propa}
 Suppose that $a$ is not identically zero in  $B_\delta^\prime  (0)$  and  $a(0)=\grad a(0)=0$. Suppose
 that $(p_n)_n$ is a sequence in $\R^N$ such that, for all $n$ and for all $x\in B_\delta^\prime  (0)$,
$$|p_n\cdot (x^\prime, a(x^\prime))|\leq C$$
for some constant $C$. Then $(p_n)_n$ is a bounded sequence. 
          \end{prop}
{\em Proof of Lemma \ref{lem1a}. } We assume first that $h=0$.  

  Note that if  $B_1(x) \cap \{ y_N = a(y^\prime)\} = \emptyset$ then  
it is sufficient to use the result of \cite{IS}. 
 So we now assume that  $B_1(x) \cap \{ y_N = a(y^\prime)\} \neq \emptyset$. 
 
Let $\rho$ and $q=q_{x,\rho}$ be so that any $u$, solution of 
$$\left\{\begin{array}{lc}
   F(D^2 u) + h(y)\cdot \nabla u = 0& {\rm in}\ B_1(x) \cap \{ y_N > a(y^\prime)\}\\
    u= \varphi& {\rm on}\  B_1(x) \cap \{ y_N = a(y^\prime)\},
    \end{array}\right.$$
satisfies $osc_{B_\rho(x)} (u(y)-q \cdot y) \leq  {\rho\over 2}$.

We argue by contradiction and suppose that for all $n$,  there exist  $x_n\in \overline{B}$  and $p_n\in \R^N$,  $|f_n|_{L^\infty (\Omega )}\leq {1\over n} $ and  a function $u_n$ with $osc   (u_n) \leq 1$ solution of 
\begin{equation}\label{hhn}
\left\{ \begin{array}{lc}
 |p_n+\nabla u_n|^\alpha F(D^2 u_n) = f_n & {\rm in}\  B\cap \{ y_N > a(y^\prime) \}\\
 u_n(y)+ p_n\cdot y =\varphi(y^\prime)  & {\rm on} \  B\cap \{ y_N=a (y^\prime)\}
 \end{array}\right.
 \end{equation}
such that  
\begin{equation}\label{bb}
{\rm osc}_{B_\rho(x_n)} (u_n(y)-q \cdot  y) > {\rho\over 2}.
\end{equation}
 
 Extract from $(x_n)_n$ a subsequence which converges to $ x_\infty \in \overline{B}\cap \{ y_N \geq a(y^\prime)\}$. 
 
  We denote in the sequel $B_\infty= \cap_{n\geq N_1}  B_1(x_n)$ which contains $B_\rho(x_\infty)$ as soon as $N_1$ is large enough.

  {\bf The case where $T$ is not straight. }
Observing that $u_n-u_n(x_n)$ satisfies the same equation as 
$u_n$,  has oscillation $1$ and is bounded,  we can  suppose that the sequence $(u_n)$ is bounded.
This, together with the boundary condition and Proposition  \ref{propa}, implies that $(p_n) $ is bounded.

So,  up to subsequences,
$u_n$  converges to some $u_\infty$,  uniformly on $\overline{B_r}\cap \{ y_N \geq a(y^\prime)\}$ for every 
$B_r\subset \subset B_\infty$, (due to Corollary \ref{compac1} ), 
and  $p_n$ converges to $p_\infty$.
Furthermore, 
 $(u_\infty, p_\infty)$ solves 
     $$\left\{ \begin{array}{lc}
 |p_\infty+\nabla u_\infty|^\alpha F(D^2 u_\infty) = 0& \ {\rm in}\ B_\infty \cap \{ y_N > a(y^\prime)\} \\
 u_\infty+ p_\infty\cdot y =\varphi(y) & {\rm on} \  B_\infty \cap \{ y_N=a(y^\prime)\}.
 \end{array}\right.$$ 
 Using Lemma 6 in \cite{IS} one gets that 
      $$\left\{ \begin{array}{lc}
F(D^2 u_\infty) = 0 & \ {\rm in}\   B_\infty  \cap \{ y_N >a(y^\prime)\}\\
u_\infty+ p_\infty\cdot y =\varphi (y) & {\rm on} \  B_\infty\cap \{ y_N=a(y^\prime)\}.
 \end{array}\right.$$ 
On the other hand,   by passing to the limit  in (\ref{bb}),
one obtains 
$$\mbox{osc}_{B_\rho( x_\infty )\cap \{ y_N >0\}} (u_\infty(y)-q \cdot y  ) \geq  {\rho\over 2}.$$  
This is a contradiction  with the assumption on $\rho$ and $q$ and it ends the proof.

\noindent {\bf The case where $T = \{ y_N = 0\}$.}
Let $p_n=p_n^\prime+p_n^Ne_N$. The boundedness of $u_n$ and the boundary condition
imply that $(p_n^\prime)_n$ is bounded. If $(p_n^N)_n$ is bounded just proceed as above. So we suppose that $p_n^N$ is unbounded.   
 
Dividing (\ref{hhn}) by $|p_n^N|^\alpha $ it becomes 
$$\left\{ \begin{array}{lc}
|{p_n \over p_n^N}+{1\over p_n^N}\nabla u_n|^\alpha F(D^2 u_n) = {f_n\over |p_n^N|^\alpha } & \ {\rm in}\   B_\infty\cap \{ y_N >0\}\\
u_n+ p_n^\prime \cdot y =\varphi (y^\prime)& {\rm on} \   B_\infty \cap \{ y_N=0\}
\end{array}\right.$$ 
Denoting   by $p^\prime $  the limit of a subsequence of $p_n^\prime$,   and  $u_\infty $ the limit of a subsequence of $( u_n)$, one gets   by passing to the limit and using Corollary \ref{compact2} 
  that 
     $$\left\{ \begin{array}{lc}
F(D^2 u_\infty) = 0& \ {\rm in}\  B_\infty \cap \{ y_N >0\}  \\
 u_\infty+ p^\prime \cdot y =\varphi(y^\prime)  & {\rm on} \ B_\infty \cap \{ y_N=0\}.
 \end{array}\right.$$ 
Passing to the limit in  (\ref{bb}) one gets that $osc_{B_\rho( x_\infty) } (u_\infty-q \cdot  y) > {\rho\over 2}$, a contradiction.
This ends the proof for $h=0$.

\bigskip
 We briefly  point out the differences in the case $h\neq 0$. 
It is sufficient to treat the case  of a  straight boundary, the other cases being as before.  
Indeed,  we already know that if the boundary is not locally straight, $(\vec p_n)_n$ is bounded. 

In the case where the boundary is locally straight, say 
$T= \{y_N=0\}$ so that the only possibly unbounded component is $p_n\cdot e_N$.

The Claim in the interior flatness Lemma \ref{flat} implies that for any open set
$D\subset B_r(x)\cap \{ y_N > a(y^\prime)\}$, $h(y)\cdot e_N=0$ for any $y\in D$.
By the arbitrariness of $D$, this implies that $h(y)\cdot e_N=0$ in $B_r(x)\cap \{ y_N > a(y^\prime)\}$.
We can now apply   Proposition \ref{liplargep} and end as in the case $h=0$. 

\bigskip

\noindent We now want to give  the 

\noindent{\em  proof of  Proposition \ref{propa}}.

\bigskip        
Since $a$ is not identically $0$,  and $a(0)=0$, $\grad a(0)=0$, there exists $i$ such that  $\partial_i a$ is 
not identically zero in $B_\delta^\prime (0)$ for some $\delta>0$.

Hence, there exist $y_1^i  \in B_\delta^\prime (0)$ and $y_2^i \in B_\delta^\prime (0)$ such that 
       $$ y_1^i  e_i+ a( y_1^i  e_i)e_N\ \mbox{ and}\ y_2^i  e_i+ a(y_2^i  e_i)e_N$$ 
are linearly independent.
Indeed, assume by contradiction that, for $y_o\neq 0$ fixed  in $B_\delta^\prime (0)$,  for all $y\in B_\delta^\prime (0)$,  there 
exists $t\in \R$ such that $(ye_i, a(ye_i)) =  t (y_o e_i, a(y_oe_i))$.
Then, 
$a(ye_i) = {a (y_o e_i) y\over y_o}$,  which implies that $a_{,i} = cte$ in $B_\delta^\prime (0)$ which is a 
contradiction.

Defining 
$E_j^i= y_j ^i e_i + a( y_j ^i) e_N$,  for $i\in [1, N-1]$ and $j = 1,2$ one gets 
$$|(p_n^i , p_n^N) \cdot (E_j^i)|\leq C\ \mbox{for}\ j=1,2.$$
 Since $E_j^i$ are two linearly independent vectors,  it
implies that $(|p_n^i | + |p_n^N|)_n$ is bounded.  
Suppose now that $a_{,i} \equiv 0$, then replacing $E_j^i$ by $e_i$ in the previous reasoning, one obtains that $|p_n^ i|$ is bounded. 
 Since $a$ is not identically zero, there exists at least one $i$ for which the first situation occurs,  so we get that for all $j$,  $(p_n^j)_n$ is bounded.

    \end{document}